\documentclass[10pt]{amsart}

\newtheorem{theorem}{Theorem}[section]
\newtheorem{lemma}[theorem]{Lemma}

\theoremstyle{definition}

\theoremstyle{remark}

\numberwithin{equation}{section}

\title[{Trace formulas for Schr\"odinger operators}]
{Trace formulas for Schr\"odinger operators\\ -- from the view point of complex analysis}

\date{\today}
\author[Hiroshi Isozaki]{Hiroshi Isozaki}
\address{Institute of Mathematics,
University of Tsukuba,  Tsukuba, 305-8571, JAPAN \\
 isozakih@math.tsukuba.ac.jp}
\author[Evgeny Korotyaev]{Evgeny L. Korotyaev}
\address{ Saint-Petersburg State
University, Russia.
 \ korotyaev@gmail.com}


\begin{document}
\maketitle

\begin{abstract}
We consider the Schr{\"o}dinger operator $-\Delta +V(x)$ in $L^2({\bf R}^3)$
with a real short-range (integrable) potential $V$. Using the associated Fredholm determinant, we present
new trace formulas, in particular, the ones in terms of resonances and eigenvalues only. We also derive expressions of the Dirichlet integral, and the scattering phase.  The proof is based on the change of view points for the above mentioned problems from the
operator theory to the complex analytic (entire) function
theory.
\end{abstract}


\section {Introduction and main results}
\subsection{Modified Fredholm determinant}

Let ${\bf B}$, ${\bf B}_1$, ${\bf B}_2$ be the set of all bounded operators, trace class and Hilbert-Schmidt class operators on $L^2({\bf R}^3)$, respectively. The norms of ${\bf B}$, ${\bf B}_1$, ${\bf B}_2$ are denoted by $\|\cdot\|, \ \|\cdot\|_1, \ \|\cdot\|_2$, respectively. For $A \in {\bf B}$, $\sigma(A)$ denotes the spectrum of $A$.

 We study the
Schr{\"o}dinger operator $H$  in $L^2({\bf R}^3)$ given by
$$
H=H_0+V,\quad  H_0=-\Delta.
$$
We put $R_0(z) = (H_0 -z)^{-1}$,
 and
\begin{equation}
Q_0(k) = |V|^{1/2}\, R_0(k^2)\, \widetilde{V}^{1/2}, \quad k \in {\bf C}_+ = \{k \in {\bf C}\, ; \, {\rm Im}\,k > 0\},
\label{S1Q0(k)}
\end{equation}
where ${\widetilde V}^{1/2}(x) = V(x)/|V(x)|^{1/2}$ for $V(x) \neq 0$, and ${\widetilde V}^{1/2}(x) = 0$ for
$V(x) = 0$. We basically assume the following condition.

\medskip
\noindent
$(C)$ : {\it The potential $V(x)$
is a real-valued $C^2$-function and satisfies
\begin{equation}
\sum_{|\alpha|\leq 2}|\partial_x^{\alpha}V(x)| \leq C(1 + |x|)^{-3-\epsilon},
\label{S1ConditonC}
\end{equation}
for some constants $C, \epsilon > 0$. The strong limit
\begin{equation}
Q_0(0) = \lim_{{\bf C}_+ \ni k\to 0}Q_0(k)
\label{S1Q(0)}
\end{equation}
exists in ${\bf B}$, moreover $- 1 \not\in \sigma(Q_0(0))$.}

\medskip
Under this condition,  the operator $H$ has absolutely continuous
spectrum $[0,\infty)$ without embedded eigenvalues and a finite number of bound states (counted with multiplicity)
$-\lambda_1 < - \lambda_2 \leq \cdots \leq - \lambda_N < 0$. We put
\begin{equation}
\sqrt{\sigma_d(H)} = \{i\sqrt{\lambda_1}, \cdots, i \sqrt{\lambda_N}\}\subset {\bf C}_+.
\label{S1Sqrtsigmad}
\end{equation}

Recall that $R_0(z)$ has the integral kernel
\begin{equation}
R_0(x-y;z) = \frac{e^{i{\sqrt z}|x-y|}}{4\pi|x-y|}, \quad \sqrt{z} \in {\bf C}_+.
\label{S1R0zkernal}
\end{equation}
Let
$R(z) = (H - z)^{-1}$ for $z \in {\bf C}\setminus{\bf R}$.
Then it is well-known that under the condition $(C)$,  for $k \in {\bf R}$, there exists the strong limit $Q_0(k) = \lim_{\epsilon\to 0}Q_0(k + i\epsilon) \in {\bf B}$. Moreover for $k \in {\bf C}_+\setminus{\sqrt{\sigma_d(H)}}$
\begin{equation}
 |V|^{1/2} R_0(k^2), \  Q_0(k)\in {\bf B}_2,\quad
  R(k^2)-R_0(k^2)\in {\bf B}_1.
\label{S1Rk2-R0k2inB1}
\end{equation}
However,
$VR_0(k^2) \not\in {\bf B}_1$ for any $k \in {\bf C}_+\setminus{\sqrt{\sigma_d(H)}}$, and we need some modification to define the determinant $\det (I+VR_0(k))$.

We introduce a notation.
For a Banach space $\mathcal X$ and a domain $D \subset {\bf C}$, let $\mathcal A(D\, ; \,\mathcal X)$ be the set of all $\mathcal X$-valued analytic functions on $D$, and
\begin{equation}
\mathcal H(D;\mathcal X) = \mathcal A(D\,;\,{\mathcal X})\cap
C(\overline{D}\, ; \, \mathcal X),
\label{S1DefineH(C+X)}
\end{equation}
where $C({\overline D}\, ; \, \mathcal X)$ is the set of all $\mathcal X$-valued continuous functions on $\overline{D}$, the closure of $D$.
 Since $Q_0\in \mathcal H({\bf C}_+\,;\,{\bf B}_2)$, we can define the modified Fredholm determinant
\begin{equation}
 D(k)=\det \left[(I+ Q_0(k))e^{- Q_0(k)}\right], \quad k \in \overline{{\bf C}_+},
\label{S1ModifiedDet}
\end{equation}
which also belongs to $\mathcal H({\bf C}_+\,;\,{\bf B}_2)$. We then have
$$
\{k \in \overline{{\bf C}_+} \, ; \, D(k) = 0 \} = \sqrt{\sigma_d(H)},
$$
which follows from the condition $(C)$ and the non-existence of positive eigenvalues of $H$.
 The function $D(k)$ is a
regularization of $\det (I+VR_0(k))$. In
fact, after removing the zeros by a suitable Blaschke product,
$D(k)$ will correspond to an element in the Hardy class, and
within this class, $D(k)$ will turn out to be a unique regularized
determinant. Let us prepare some notation to make it precise.

Since $\|Q_0(k)\|_{2}=o(1)$ as ${\rm Im}\, k \to \infty$,
 we can define the branch
\begin{equation}
\log D(k) =  \rho(k) + i\phi(k),
\label{S1logDkrhophi}
\end{equation}
\begin{equation}
 \rho(k) =  \log|D(k)|, \quad \phi(k) = {\rm arg}\, D(k),
\label{S1rhokphik}
\end{equation}
so that $\log D(k)=o(1)$ as ${\rm Im}\, k\to \infty$.

Our aim is to deduce the following terms from $\log D(k)$:
\begin{equation}
 \alpha_{-1}=\frac{1}{4\pi}\int_{{\bf R}^3}V(x)dx,
\label{S1alpha-1}
\end{equation}
\begin{equation}
 \alpha_0=\frac{1}{(4\pi)^2}\int_{{\bf R}^3} V^2(x)dx,
\label{S1alpha0}
\end{equation}
\begin{equation}
\alpha_1=\frac{1}{3(4\pi)^3}\int_{{\bf R}^3}
\left(|\nabla V(x)|^2+2V^3(x)\right)dx.
\label{S1alpha1}
\end{equation}

Let us remove the zeros from $D(k)$. We define
\begin{equation}
B(k)=\prod_{j=1}^N\left(\frac{k-i\sqrt{\lambda_j}}{ k+i\sqrt{\lambda_j}}\right),
\label{S1B(k)}
\end{equation}
and put
\begin{equation}
D_B(k) = \frac{D(k)}{B(k)}.
\label{S1DB(k)}
\end{equation}
We then have
\begin{equation}
\log D_B(k) = \log D(k) - \log B(k),
\label{S1logDB(k)}
\end{equation}
where $\log B(i\tau)=o(1)$ as $\tau \to
\infty$. The function $\log B(k)$ is analytic in the domain ${\bf C}\setminus
\{i\tau\, ; \, \tau \in [-\sqrt{\lambda_1},\sqrt{\lambda_1}]\}$ and has the following expansion
\begin{equation}
 i\log B(k)=\frac{\beta_0}{k}-\frac{\beta_2}{k^3}+\frac{\beta_4}{k^5} + \cdots,  \quad {\rm
if} \quad|k|>\sqrt{\lambda_1}>0,
\label{S1logB(k)}
\end{equation}
 where
\begin{equation}
\beta_n = \frac{2}{n+1}\sum_{j=1}^N
(\sqrt{\lambda_j})^{n+1},\quad  n\neq -1.
\label{S1betan}
\end{equation}
It is convenient to put
\begin{equation}
\gamma_{n}=\alpha_{n}-(-1)^n\beta_{2n}, \quad n=-1, 0,1,\cdots.
\label{S1gamman}
\end{equation}
Note that we have $i\log B(0)=\pi N$ and $\log D_B(k) \in \mathcal H({\bf C}_+;{\bf C})$.


\subsection{Scattering phase and modified determinant}
The condition (C) implies the existence and completeness of wave operators for the pair $H_0$ and $H$. The
associated  scattering matrix  $S(k)$ has the property $S(k)-I\in{\bf B}_1$ for any $k\in {\bf R}$. Moreover, there exists a unique function
(the scattering phase) $\phi_{sc}(k)$, which is odd on ${\bf R}$,
continuous on ${\bf R}_+$ and satisfies $\phi_{sc}(+0)=-\pi N$. This
function $\phi_{sc}(k)$ is connected with the scattering matrix $S(k)$
by the the Birman-Krein identity (see p. 6 of  \cite{BY})
\begin{equation}
 \det S(k)=e^{- 2i \phi_{sc}(k)}, \quad {\rm a.e.} \quad k>0.
\label{S1detSk}
\end{equation}
We recall the following important identity from  \cite{N}:
\begin{equation}
 \det S(k)=\frac{\overline{D(k)}}{D(k)}e^{-2ik\alpha_{-1}},
\label{S1detSkandDk}
\end{equation}
which yields
$$
\phi_{sc}(k)=k\alpha_{-1}+ \arg D(k), \quad k>0.
$$


\subsection{Main results}
Now we can state our first main result on the trace formula.


\begin{theorem}
(1) We have the following asymptotic expansion as $|k| \to \infty$
\begin{equation}
- i\log D(k) = - \frac{\alpha_0}{k} - \frac{\alpha_1}{k^3} + O\Big(\frac{1}{k^4}\Big),
\label{S1ExpandlogDk}
\end{equation}
\begin{equation}
- i\log D_B(k) = - \frac{\gamma_0}{k} - \frac{\gamma_1}{k^3} + O\Big(\frac{1}{k^4}\Big),
\label{S1ExpandlogDBk}
\end{equation}
both being uniform with respect to ${\arg}\,k \in [0,\pi]$. \\
 \noindent
 (2) The following trace formulas hold:
\begin{equation}
\frac{1}{\pi}\int_{{\bf R}}t\left[{\rm arg}\,D(t) + \frac{\alpha_0}{t}\right]dt = \beta_1,
\label{S1Traceformula1}
\end{equation}
\begin{equation}
\frac{1}{\pi}\int_{{\bf R}}t^3\left[{\rm arg}\,D(t) + \frac{\alpha_0}{t} + \frac{\alpha_1}{t^3}\right]dt = - \beta_3,
\label{S1Traceformula2}
\end{equation}
\begin{equation}
\frac{1}{\pi}\int_{{\bf R}}\log |D(t)|dt = - \gamma_0,
\label{S1Traceformula3}
\end{equation}
\begin{equation}
\frac{1}{\pi}\int_{{\bf R}}k^2\log |D(t)|dt = - \gamma_1,
\label{S1Traceformula4}
\end{equation}
\begin{equation}
\frac{1}{\pi}\int_{{\bf R}}\frac{\log |D(t)|- \log|D(0)|}{t^2}dt
+ \phi'_{sc}(+0) =  \gamma_{-1}.
\label{S1Traceformula5}
\end{equation}
\end{theorem}

\medskip
\noindent
{\bf Remark.} (i) The formulas (\ref{S1Traceformula1}), (\ref{S1Traceformula2}) were proved by Buslaev
\cite{B}. The trace formulas (\ref{S1Traceformula3}) and  (\ref{S1Traceformula4}) are new. The equality (\ref{S1Traceformula5}) is an analogue of (\ref{S1Traceformula3}), (\ref{S1Traceformula4}). In fact, it follows from the asymptotics
(\ref{S1ExpandlogDk}).

\noindent
(ii) Higher regularity of  $V$ implies more trace formulas.

\noindent
(iii) Buslaev mainly considered the phase
$\phi_{sc}$. In the present paper, in addition, the trace formulas for
the conjugate function $-\log|D(k)|$ are proved, which gives rise to a more complete
result.

\medskip
The scattering phase
$\phi_{sc}$ defined by (\ref{S1detSk}) is harmonic in ${\bf C}_+\setminus\{i\tau\, ; \, \tau \in [0,\sqrt{\lambda_1}]\}$. Thanks to the analyticity of $D(k)$ in ${\bf C}_+\setminus\{i\tau\, ; \, \tau \in [0,\sqrt{\lambda_1}]\}$, we can construct the conjugate function  $-\log|D(k)|$ in the same region. In fact we have 2
alternative definitions of the conjugate harmonic function $-\log|D(k)|$: either

\smallskip
(1) First define $\phi_{sc}$ directly from the $S$-matrix and then  determine $-\log|D_B(k)|$
from ${\rm arg}\,D_B(k)$ using the Hilbert transform, or

\smallskip
 (2) Define implicitly $-\log|D(k)|$ in terms of the determinant $D(k)$.

\smallskip

 In fact, there are no works on the conjugate function $-\log|D(k)|$ of the
scattering phase  $\phi_{sc}$ in the case dimension $>1$. Only
in the 1 dimensional case there are a lot of papers devoted to the
trace formulas and the conjugate function $-\log|D(k)|$, since it plays an important role in spectral theory, inverse
problems, non-linear equations, etc. For example, $-\log|D(k)|$ is the action variable
for the KDV equation (see \cite{FZ}). For periodic potentials such
identities were obtained in \cite{KK1}, \cite{K2}, \cite{MO} and were
used to get double sided estimates of potentials in terms of
spectral data to solve the inverse problem \cite{KK2}.

\medskip

We formulate the second result on the trace formulas.


\begin{theorem}
The Dirichlet integral of $\log D_B(k)$ has the following form
\begin{equation}
\frac{1}{\pi}\iint_{{\bf C}_+}\Big|\frac{d}{dk}\log D_B(k)\Big|^2dtd\tau +S_0= m_B\gamma_0,
\label{ddklogDBk}
\end{equation}
where, $k = t + i\tau$, and
$$
m_B =-\min_{t\in {\bf R}}\frac{d}{dt}{\rm arg}\,D_B(t), \quad
S_0=-\frac{1}{\pi}\int_{{\bf R}}\log|D(t)|(m_B+\frac{d}{dt}{\rm arg}\,D_B(t))dt.
$$
Furthermore, if
$N=0$, then $\log|D(0)| < 0$.
\end{theorem}

\medskip
\noindent
{\bf Remark.} (i) The functions $\log D_B(k),
k^{-1}{\rm arg}\,D_B(k)$ are continuous on ${\bf R}$, since ${\rm arg}\,D_B(k)$ is even on the real line and ${\rm arg}\,D_B(0)=0$. Moreover, by (\ref{S1ExpandlogDBk}), they belong to
$L^2({\bf R})$. By the property of the Hilbert transform we have
$$
\frac{1}{\pi}\int_{{\bf R}}\frac{{\rm arg}\,D_B(t)}{t}dt = \log|D(0)|.
$$
If $V \geq 0$, then $\phi_{sc}(t) > 0$ for any $t > 0$. However, ${\rm arg}\,D_B(t) < 0$  for some $t > 0$.

\noindent
(ii) Let $\lambda(k) = m_Bk - i\log D_B(k)$, and consider the curve $\Gamma = \{\lambda(t)\, ; \, t \in {\bf R}\}$. Then $S_0$ is the area of the region surrounded by $\Gamma$ and ${\bf R}$.

\noindent
iv) Note that for the Hill operator there are identities of the type
(\ref{ddklogDBk}) proved in \cite{KK1}, \cite{K2}. Roughly speaking in
the case of the Hill operator we have $S_0=0$, while for the Schr\"odinger
operator $-\Delta +V(x)$ in $L^2({\bf R}^3)$ with a real short-range
potential $V$, the additional term $S_0$ appears.

\medskip
We turn to the resonance. If the function $D(k)$ is entire, then
it has $N\geq 0$ zeros $i\sqrt{\lambda_1},\cdots,i \sqrt{\lambda_N}$ in the upper half-plane, and (in general) an infinite
number of zeros in the lower half-plane ${\bf C}_- = \{k \in {\bf C}\, ; \, {\rm Im}\,k < 0\}$, which we can arrange as $0< |k_{N+1}|\leq |k_{N+2}|\leq \cdots$. It is convenient to put
$$
k_n = i\sqrt{\lambda_n}, \quad 1 \leq n \leq N.
$$
 We call
the zeros of $D(k)$ in ${\bf C}_-$ {\it resonances}  of $H$.
Let
\begin{equation}
Q(k) = |V|^{1/2} R(k^2)\widetilde{V}^{1/2} .
\label{S1Q(k)}
\end{equation}


\begin{theorem}
(1) The function $B^2(k)\det S(k), k\in {\bf R}$, has an analytic
continuation onto whole of ${\bf C}_+$ if and only if the function $D(k)$ is entire.

 We assume that  $D(k)$ is entire.

\noindent
(2) $\det S(k), k\in {\bf R}$,  has a meromorphic continuation
 onto whole ${\bf C}$ and the following formula holds
\begin{equation}
e^{2ik\alpha_{-1}}B^2(k)\det S(k)=\frac{D_B(-k)}{D_B(k)}, \quad
k\in {\bf C}\setminus \{0, k_1, k_2,\cdots\}.
\label{S1D0-kD0k}
\end{equation}
(3)  Suppose  that the function  $D(k)$ has a finite  order and
a finite number of zeros. Then $V = 0$.

\noindent
(4) Suppose $z^2 \not\in \sigma_d(H)$. Then
$z$ is a pole of $\det S(k)$ (of
multiplicity  $m\geq 1$) if and only if  $z$ is a zero of $D(k)$ (of the same
multiplicity $m\geq 1$).

\noindent
(5) Suppose as a $\bf B_2$-valued function, $Q_0(k)$ has an analytic
continuation to the whole plane ${\bf C}$. Then $Q(k)\, : \, {\bf C}_+\setminus \sqrt{\sigma_d(H)} \to {\bf B}$
has a meromorphic continuation to whole ${\bf C}$. Moreover, $z$ is a zero
of $D(k)$, if and only if  $z$ is a pole of $Q(k)$ with the same
multiplicity.
\end{theorem}

\medskip
\noindent
{\bf Remark.} It is well known that if $V(x)$ is compactly supported, more generally if
$\sup_{x\in {\bf R}^3}|V(x)|e^{c|x|}<\infty$ for any $c > 0$, then as a ${\bf B}_2$-valued function, $Q_0(k)$
has an analytic continuation into the whole plane ${\bf C}$, and  $D(k)$ is
entire.


\begin{theorem}
Assume that ${\rm supp}\, V \subset \{|x|\leq 1\}$. Then
 $D(k)$ is an entire function of order $\leq 3$, and
\begin{equation}
 |D_B(k)|\leq \sup_{t\in{\bf R}}|D(t)|, \quad k \in \overline{{\bf C}_+},
\label{S1DBKinfhk}
\end{equation}
\begin{equation}
 |D(k)|\leq C_0e^{C_0|k|^3}, \quad k\in {\bf C}_-,
\label{S1DkeC0k3}
\end{equation}
for some constant $C_0$. Moreover
\begin{equation}
 D(k)=D(0)e^{{\mathcal P}(k)}\lim_{r\to +\infty}\prod_{|k_n|\leq
r}\left(1-\frac{k}{k_n}\right) e^{\frac{k}{k_n}+\frac{k^2}{2 k_n^2}+
\frac{k^3}{3 k_n^3}}, \quad k\in {\bf C},
\label{S1Dkinfiniteprodut}
\end{equation}
uniformly on any compact subset of ${\bf C}$,
where
\begin{equation}
{\mathcal P}(k)=i\frac{\phi'''(0)}{6}k^3+\frac{h''(0)}{2}k^2+i\phi'(0)k,
\label{S1mathcalPk}
\end{equation}
with $\phi(k)$ defined in (\ref{S1rhokphik}),  and
\begin{equation}
 \frac{D'(k)}{D(k)}={\mathcal P}'(k)+\sum_{n\geq1} \frac{k^3}{
k_n^3(k-k_n)},
\label{S1D7koverD(k)}
\end{equation}
the convergence being uniform on any compact subset of ${\bf C}\setminus \{0, k_1,k_2,\cdots\}$.
\end{theorem}

\medskip
\noindent
{\bf Remark.} (i) In \cite{K1} the corresponding 1-dimensional case is discussed.

\noindent
(ii) The 1-dimensonal counter part of $D(k)$ is of exponential  type (see \cite{Z1}), and we can make use of the   theory  of functions of exponential
type (see \cite{K1}, where the results in \cite{Koo} was used).  In dimension 3 such nice properties are unknown, and we stress {\it the fundamenal importance of the problem of determining the order of $D(k)$}.

\medskip

Recall Krein's trace formula for the pair of operators $H_0, H$ and
$f \in C_0^{\infty}(\bf R)$:
\begin{equation}
  {\rm Tr}\,(f(H)-f(H_0))=\int_{\bf R}\xi (E) f'(E)dE,
\label{S1Kreintrace}
\end{equation}
where the spectral shift function $\xi $ has the following form (the
Birman-Krein formula)
\begin{eqnarray}
\xi(E) &=&
\left\{
\begin{split}
& \frac{1}{\pi}\phi_{sc}(\sqrt E) \quad {\rm if} \quad E>0, \\
&  - \int _{-\infty}^E \sum _{n=1}^N\delta (t+\lambda_n)dt \quad {\rm if} \quad  E<0,
\end{split}
 \right.\\
\xi(+0) &=& -N.
\end{eqnarray}


\begin{theorem}
Let $V$, $\phi(k)$ and $\mathcal P(k)$ be as in Theorem 1.4.
Then for any $f\in C_0^{\infty}(\bf R )$ the following identities hold:
\begin{equation}
{\rm Tr}\,\left(f(H)-f(H_0)\right)=\sum_{j=1}^{N}f(k_j^2)-\frac{\alpha_{-1}}{\pi}\int_{0}^{\infty}f(t^2)dt
-\frac{1}{\pi}\int_{0}^{\infty}f(t^2)\phi'(t)dt,
\label{S1TraceformulainTheorem15}
\end{equation}
\begin{equation}
\phi'(t)={\rm Im}\,{\mathcal P}'(t)+{\rm Im}
\sum_{n\geq1} \frac{t^3}{k_n^3(t-k_n)}, \quad t \geq 0,
\label{S1ph7kinTheorem15}
\end{equation}
where  the series converges absolutely and uniformly on any compact set
of $[0,\infty)$.
\end{theorem}

\medskip
\noindent
{\bf Remark.} The equality (\ref{S1ph7kinTheorem15})  is the Breit-Wiegner formula for the
resonance (see p. 53 of \cite{RS2}).

\medskip

The importance of trace formula in non-linear equations,
 inverse problems, spectral theory has already been discussed in many papers,  see \cite{DK}, \cite{FZ},
 \cite{KK1}, \cite{LW}
 and references therein. The trace
formula associated with the scattering phase (the spectral shift
function) $\phi_{sc}$ was derived by Buslaev \cite{B} and later
it is repeatedly studied by various authors, see \cite{C}, \cite{N},
\cite{O}, \cite{R} etc. For example, using the modified
determinant,  Newton \cite{N} gave a new proof of the Levinson Theorem.
The trace formulas for Stark operators were obtained in \cite{KP1} and for the case 2D magnetic fields see
\cite{KP2}.

The resonance is also a broadly discussed subject (see
\cite{Me}, \cite{Z2}, \cite{Fr}, \cite{SZ} etc). Many related problems which have applications to physics are still open.  We think that some basic
problems for the distribution of resonances have connections to open problems in the theory of entire
functions. The works cited in \cite{Koo} will serve as a first step
in this direction. For example, the Levinson Theorem (about zeros of
entire functions, see \cite{Koo} , \cite{L}) gives only the first term in the
asymptotics of the number of resonances in dimension 1. We expect that the
next term in these asymptotic estimates could be determined  by
applying the theories of entire functions and conformal mapping.


\section {Trace formulas}

Let us recall some well-known facts. Let $A, B\in {\bf B}$ and $AB, BA\in
{\bf B}_1$. Then
\begin{equation}
{\rm Tr}\, AB={\rm Tr}\, BA,
\label{S2TrAB=TrBA}
\end{equation}
\begin{equation}
\det (I+ AB)=\det (I+BA).
\label{S2det1+AB}
\end{equation}
Suppose for a domain $D \subset {\bf C}$, $\Omega(z) \in {\mathcal A}(D\, ;\, {\bf B}_1)$ satisfies $-1 \not\in \sigma(\Omega(z))$ for any $z\in
D$. Then for  $F(z)=\det (I+\Omega (z))$ we have
\begin{equation}
 F'(z)= F(z){\rm Tr}\,(I+\Omega (z))^{-1}\Omega '(z).
\label{S2F'z}
\end{equation}
We need the following fact for the trace class operators
(see Theorem XI.21 \cite{RS1}):
\begin{equation}
f(x)g(-i\nabla )\in {\bf B}_1,\quad f, g\in L^{2,\delta}({\bf R}^n),
\quad \delta > n/2,
\label{S2fxgidela}
\end{equation}
where $f\in L^{2,\delta}({\bf R}^n)$ means that
$\int_{{\bf R}^n}(1 + |x|)^{2\delta}|f(x)|^2dx < \infty$.
Recall that $\sqrt{\sigma_d(H)} =\{i\sqrt{\lambda_1},\cdots,i\sqrt{\lambda_N}\}$, and $Q_0(k)$ and $Q(k)$ defined in (\ref{S1Q0(k)}) and (\ref{S1Q(k)}).
We use the following formula
\begin{equation}
 (I+Q_0(k))(I-Q(k))=I,\quad
 k\in {\bf C}_+\setminus \sqrt{\sigma_d(H)}.
\label{S2I+Q0I+Q}
\end{equation}
We put
\begin{equation}
Q_B(k) = B(k)Q(k).
\label{S2QB(k)}
\end{equation}
Recall also that for a Banach space $\mathcal X$, $\mathcal H({\bf C}_+\, ; \, \mathcal X)$ is the set of all $\mathcal X$-valued continuous functions on $\overline{{\bf C}_+} = \{k \in {\bf C}\, ; \, {\rm Im}\;k \geq 0\}$, which are analytic on ${\bf C}_+$.


\begin{lemma}
The operator-valued functions $Q_0, Q_0', Q_B, Q_B'$ belong to $\mathcal H({\bf C}_+\,;\,{\bf B}_2)$. Moreover, we have:
\begin{equation}
\sup_{k\in {\bf C}_+}\left(\|Q_0'(k)\|_{2}+\|Q_0(k)\|_{2}\right)<\infty,
\label{S2Q0Q0'estimate}
\end{equation}
\begin{equation}
 \sup_{k\in {\bf C}_+}\,|k|\|Q_0(k)\|<\infty,
\label{S2kQ0kestimate}
\end{equation}
\begin{equation}
VR_0(k)^2,\ \ Q_0'(k)\in {\bf B}_1, \quad  k\in {\bf C}_+,
\label{S2VR0kQ0'k}
\end{equation}
\begin{equation}
{\rm Tr}\ Q_0'(k)=i\alpha_{-1}, \quad k\in {\bf C}_+.
\label{S2TraceQ0'k}
\end{equation}
If $V$ satisfy (\ref{S1ConditonC}) with $\epsilon >1$, then
$Q_0''(k) \in {\mathcal H}({\bf C}_+\, ; \, {\bf B}_2)$.
\end{lemma}

Proof. It is well known that $Q_0, Q_0', Q_B \in {\mathcal H}({\bf C}_+\, ; \, {\bf B}_2)$ and
(\ref{S2Q0Q0'estimate}) holds. Note that it is a simple fact and it follows from (\ref{S1R0zkernal}).
Using (\ref{S2I+Q0I+Q}), we have
$$
 Q'(k)= 2k|V|^{1/2} R(k)^2{\widetilde V}^{1/2}= (I-Q(k))Q_0'(k)(I-Q(k)),
$$
which yields $Q_B'\in {\mathcal H}({\bf C}_+\,;\, {\bf B}_2)$. The estimate (\ref{S2kQ0kestimate}) is proved in
\cite{GM}, and (\ref{S2fxgidela}), (\ref{S1Rk2-R0k2inB1}) imply (\ref{S2VR0kQ0'k}).
Using (\ref{S2VR0kQ0'k}), we obtain
$$
{\rm Tr}\, Q_0'(k)= 2k{\rm Tr}\,(VR(k)^2)=\frac{i}{4\pi }\int_{{\bf R}^3} V(x)dx=i\alpha_{-1},\quad
 k\in {\bf C}_+.
$$
If $\epsilon>1$, then (\ref{S1Rk2-R0k2inB1}) shows $Q_0''\in {\mathcal H}({\bf C}_+\,;\,{\bf B}_2)$. \qed

\medskip
 The function $D(k)$ is real on $\{i\tau \, ; \, \tau \in (0,\infty)\setminus[0,\sqrt{\lambda_1}]\}$, and we have
\begin{equation}
 \overline{D(k)}=D(-\overline{k}), \quad  \phi(k)=- \phi(-\overline{k}), \quad
\rho(k)=\rho(-\overline{k}),
\label{S2batDkDbark}
\end{equation}
for $k\in \overline{\bf C_+}\setminus\{i\tau\, ; \, \tau \in [0,\sqrt{\lambda_1}]\}$.
The logarithmic derivetive of $D(k)$ has the following form (see \cite{GK}),
\begin{equation}
\frac{d}{dk}\log D(k) = - {\rm Tr} \left[ Q(k)Q_0'(k)\right],
\quad k\in {\bf C}_+\setminus\{i\tau\, ; \, \tau \in (0,\sqrt{\lambda_1}]\}.
\label{S2logDerDk}
\end{equation}
In fact, this follows from (\ref{S2TrAB=TrBA})-(\ref{S2F'z}), and (\ref{S2I+Q0I+Q}).


\begin{lemma}
The functions $\log D_B(k)$ and $\dfrac{d}{dk}\log D_B(k)$ belong to $\mathcal H({\bf C}_+\,;\,{\bf C})$,
and the following identity and the estimate hold:
\begin{equation}
C_0 := 2\sup_{k\in{\bf C}_+} |k|\|Q_0(k)\| < \infty,
\label{S2supkQ0k}
\end{equation}
\begin{equation}
-\log D(k) = \sum _{n=2}^{\infty}\frac{{\rm Tr}\,\big(-Q_0(k)\big)^n}{n},
\quad |k| > C_0, \quad k \in {\bf C}_+,
\label{S2Qkexpand}
\end{equation}
where the series converges absolutely and uniformly, and
\begin{equation}
 \Big|\log D(k)+\sum _{n=2}^{N+2}\frac{{\rm Tr}\,(-Q_0(k))^n}{n}\Big|\le
\frac{\|Q_0(k)\|^{N+1}\|Q_0(k)\|_{2}^2}{N+3},\quad |k| > C_0, \quad
 k\in {\bf C}_+.
\label{S2logDkremainder}
\end{equation}
for any $N\geq 0$. If $V$ satisfy (\ref{S1ConditonC}) with
$\epsilon >1$, then $\dfrac{d^2}{dk^2}\log D_B(k) \in \mathcal H({\bf C}_+\,;\,{\bf C})$.
\end{lemma}

Proof.   Lemma 2.1 and the formula (\ref{S2logDerDk}) imply that $\log D_B(k)$ and $\dfrac{d}{dk}\log D_B(k)$ belong to $\mathcal H({\bf C}_+\,;\,{\bf C})$.
We denote the series in (\ref{S2Qkexpand}) by $F(k)$.  Since
\begin{equation}
 |{\rm Tr}\, Q_0^n(k)|\leq \|Q_0(k)\|_{2}^2\|Q_0(k)\|^{n-2}\leq
\|Q_0(k)\|_{2}^2 \epsilon_k^{n-2}, \quad \epsilon_k=\|Q_0(k)\|<\frac{1}{2},
\label{S2TrQ0kn}
\end{equation}
 $F(k)$  converges absolutely and
uniformly, and is anlytic in $|k|> C_0$.
Moreover, differentiating (\ref{S2Qkexpand}) and using (\ref{S2I+Q0I+Q}),
 we have
$$
F'(k)=-i\sum_{n=2}^{\infty}{\rm Tr}\,(-Q_0(k))^{n-1}Q_0'(k)= i{\rm Tr}\,
Q(k)Q_0'(k),  \quad  |k|>C_0.
$$
Then we have $F(k) = - i\log D(k)$, since $F(i\tau)=o(1)$ as $\tau \to \infty$. Using (\ref{S2Qkexpand}) and (\ref{S2TrQ0kn}),  we obtain (\ref{S2logDkremainder}).

Let, in addition, $V$ satisfy (\ref{S1ConditonC}) with
$\epsilon >1$. It is
sufficient to consider $\dfrac{d}{dk}\log D(k)$ near the real line.
Using (\ref{S2logDerDk}), we have
$$
\frac{d^2}{dk^2}\log D(k)=- {\rm Tr}\,\Big(Q'(k) Q_0'(k)+ Q(k)Q_0''(k)\Big),
\quad k\in {\bf C}_+.
$$
Then by Lemma 2.1, the two terms are analytic in ${\bf C}_+$ and
continuous up to ${\bf R}$. Hence we have the last assertion of the lemma. \qed

\medskip
The following lemma will be proved in \S 4.


\begin{lemma}
We have the following asymptotic expansion as $|k| \to \infty$, $k \in {\bf C}_+$:
\begin{equation}
i{\rm Tr}\,\left(\frac{Q_0(k)^2}{2} - \frac{Q_0(k)^3}{3}\right) =
- \frac{\alpha_0}{k} - \frac{\alpha_1}{k^3} + O(\frac{1}{k^4}),
\label{S2TrQ0k2-Q0k3asymp}
\end{equation}
\begin{equation}
{\rm Tr}\,Q_0(k)^p = O(k^{-4}), \quad {\rm Tr}\, Q_0^{p-1}(k)Q_0'(k) = O(k^{-4}), \quad p = 4, 5.
\label{S2p=4,5}
\end{equation}
\end{lemma}

\medskip
Now we can compute the asymptotics of $\log D(k)$ and $\log D_B(k)$.


\begin{theorem}
We have the asymptotics (\ref{S1ExpandlogDk}), (\ref{S1ExpandlogDBk}) in Theorem 1.1, moreover
\begin{equation}
-i\frac{d}{dk}\log D(k) = \frac{\alpha_0}{k^2} + O(\frac{1}{k^4}).
\label{S2logDkprimeasymp}
\end{equation}
We also have
\begin{equation}
\inf_{t\in{\bf R}}|D(t)| \leq |D_B(k)| \leq \sup_{t\in{\bf R}}|D(t)|, \quad \forall k \in \overline{{\bf C}_+}.
\label{S2h(k)infsup}
\end{equation}
\end{theorem}

Proof. Using Lemma 2.2,  we decompose $-i\log D(k)$ as
$T_{\leq 5}(k) + T_{>5}(k)$, where
$$
 T_{\leq 5}(k)=i\,{\rm Tr}\,
\left(\frac{Q_0(k)^2}{2}-\frac{Q_0(k)^3}{3} + \frac{Q_0(k)^4}{4}-\frac{Q_0(k)^5}{5}\right),
$$
$$
T_{>5}(k)=i\sum_{n>5}\frac{{\rm Tr}\,(-Q_0(k))^n}{n}.
$$
Lemma 2.1 implies
$$
 |T_{>5}(k)|\leq \|Q_0(k)\|_{2}^2\|Q_0(k)\|^4 \leq C|k|^{-4},\quad k\in \overline{\bf C_+},\quad |k|\to\infty
$$
for some $C>0$.
The asymptotics (\ref{S1ExpandlogDk}), (\ref{S1ExpandlogDBk}) then follow
from  this estimate and Lemma 2.3. The proof of (\ref{S2logDkprimeasymp}) is similar, since we have (\ref{S2Q0Q0'estimate}).

We let
$$
\log D_B(k) = \log|D_B(k)| + i \arg D_B(k)=:\rho_B(k) + i\phi_B(k).
$$
Then $\rho_B(k)$ is a real harmonic functiuon. In view of (\ref{S1ExpandlogDBk}),
we have $\log D_B(k) = O(k^{-1})$ as $|k| \to \infty$, $k \in \overline{{\bf C}_+}$. Then the maximum principle implies
$$
\inf_{t\in{\bf R}} \rho_B(t) \leq \rho_B(k) \leq \sup_{t\in{\bf R}}\rho_B(t), \quad
k \in \overline{{\bf C}_+}.
$$
Since $\rho_B(t) = \rho(t)$ for $t \in {\bf R}$, we obtain (\ref{S2h(k)infsup}). \qed

\medskip
 For $k \in {\bf C}_+\setminus\sqrt{\sigma_d(H)}$, we put
\begin{equation}
J_0(k) = 1 + Q_0(k), \quad J(k) = 1 - Q(k),
\label{S2J0kJk}
\end{equation}
\begin{equation}
S_0(k)=J_0(-\overline{ k})J(k).
\label{SS2S0k}
\end{equation}
By (\ref{S2I+Q0I+Q}), we have
$$
S_0(k)=
1-(Q_0(k)-Q_0(-\overline{k}))(1-Q(k)).
$$
Then the operator-valued function $S_0(\cdot)-I :{\bf C}_+\to {\bf B}_1$ is
continuous up to ${\bf R}$, since the function $\det (I+\cdot)$ is
continuous in the trace norm. Then we obtain the well known formula
(see Theorem XI.42 \cite{RS1} and (\ref{S2det1+AB}))
$$
\det S_0(k)=\det S(k), \quad  k\in {\bf R}.
$$
We represent $\det S(k)$ in
terms of $D(k)$ and give an alternative proof of (\ref{S1detSkandDk}).


\begin{lemma}
 The equality (\ref{S1detSkandDk}) holds.
 \end{lemma}

Proof.  Take $z = i\tau$, $\tau \in {\bf R}\setminus\sqrt{\sigma_d(H)}$ arbitrarily, and  define the modified determinant
$$
\mathcal D (k)=\det\left[J_0(k)J(z)\right], \quad  k\in {\bf C}_+,
$$
It is well defined since $J(\cdot )J_1(z)-I\in \mathcal H({\bf C}_+\,;\,{\bf B}_1)$.
The function $\mathcal D(k)$ is analytic in ${\bf C}_+$, with $N$ zeros
(counted with multiplicity) $k_1,\cdots,k_N\in i{\bf R}_+$ and $\mathcal D(z)=1$. We put
$$
f(k)=\frac{D(k)}{D(z) }e^{{\rm Tr}\,(Q_0(k)-Q_0(z))}, \quad
k \in \overline{{\bf C}_+},
$$
 and show
 $$
 \mathcal D(k) = f(k), \quad k \in \overline{{\bf C}_+}.
 $$
Using (\ref{S2TrAB=TrBA}), (\ref{S2F'z}) and
(\ref{S2fxgidela}), we have
$$
\frac{\mathcal D'(k)}{\mathcal D(k)}= {\rm Tr}\,\left[J_0(k)J(z)\right]^{-1}Q_0'(k)J(z)
={\rm Tr}\,J(k)Q_0'(k),\quad k \in {\bf C}_+\setminus\sqrt{\sigma_d(H)}.
$$
By the similar argument, (i.e. using (\ref{S2TrAB=TrBA}), (\ref{S2F'z}), (\ref{S2fxgidela}), and (\ref{S2Q0Q0'estimate})),
 we obtain
$$
\frac{f'(k)}{f(k)}={\rm Tr}\,\left[-Q(k)Q_0'(k)+Q_0'(k)\right]=
{\rm Tr}\,J(k)Q_0'(k).
$$
Then $\mathcal D=f$, since$f$ and $\mathcal D $ satisfy the same equation and $f(z)=\mathcal D(z)=1$.

Using $J_0(k)J(k)=I$, we rewrite $\det S_0(k)$ in the form
$$
\det S_0(k)=\det \left[J_0(-\overline{ k})J(z) \right]\det
\left[J(z)^{-1}J(k)\right]= \frac{\mathcal D(-\overline{ k})}{\mathcal D(k) }=
\frac{\overline{\mathcal D(k)}}{\mathcal D(k)},
$$
for $k\in {\bf C}_+\setminus
\sqrt{\sigma_d(H)}$,
since $\mathcal D(k)$ is real on $i{\bf R}_+$.
This equality and $\det
S_0(k)=\det S(k), \ k\in {\bf R}$, yield (\ref{S1detSkandDk}), since
(\ref{S2TraceQ0'k}) gives
${\rm Tr} (Q_0(k)-Q_0(z))=i(k-z)\gamma_{-1}$. \qed

\medskip
We prepare some simple equalities. We put
\begin{equation}
B_R(k) = {\rm Re}\,B(k), \quad B_I(k) = {\rm Im}\,B(k).
\label{S2BRBI}
\end{equation}
Then we have
\begin{equation}
B_I'(0) =-{\rm Re}\, i\frac{B'(0)}{B(0)}= - {\rm Re}\, \sum_{j=1}^N\left(
\frac{i}{k-k_j}-\frac{i}{k+k_j} \right)\Big|_{k=0}=
 2\sum _{j=1}^N\frac{1}{|k_j|},
 \label{S2B'atk=0}
\end{equation}
\begin{equation}
\log D_B(0) = \log|D_B(0)|,
\label{S2logBatk=0}
\end{equation}
\begin{equation}
 -i\frac{d}{dk}\log D_B(k)\Big|_{k=0} =
\phi_B'(0) = \phi'(+0) - 2\sum_{j=1}^N\frac{1}{|k_j|},
\label{S2logDBdiffatk=0}
\end{equation}
since $\phi_B$ is odd and $\rho_B$ is even on ${\bf R}$.

\medskip

We prove the first main result about the trace formulas.

\medskip
\noindent
{\bf  Proof of Theorem 1.1}.
We put
$$
\Psi(k) = - i\log D(k), \quad \Psi_B(k) = -i\log DB(k).
$$
Due to Lemmas 2.1 and 2.2, $\Psi_B$ is continuous on $\overline{\bf C_+}$.
The asymptotics in Theorem 1.1 (1) have been proved in Theorem 2.4

In order to show the equalities (\ref{S1Traceformula3}), (\ref{S1Traceformula4}) and (\ref{S1Traceformula5}), we need the
following simple result. Assume that a function $f$ satisfies the
following condition:
\begin{equation}
 f\in \mathcal H({\bf C}_+\,;\,{\bf C}), \quad  {\rm Im}\,f(\cdot+i0)\in L^1({\bf R}), \quad
f(i\tau)=-\frac{Q_f+o(1)}{i\tau},\quad \tau\to \infty.
\label{S2Condtiononf}
\end{equation}
Then
\begin{equation}
Q_f=\frac{1}{\pi}\int_R{\rm Im}\,f(t)dt.
\label{S2Imfkint}
\end{equation}

 We apply this result to
$\Psi_B(k)$, which satisfies the
condition (\ref{S2Condtiononf}) by virtue of Theorem 1.1,  and we have (\ref{S1Traceformula3}).

Applying similar arguments to the function $k^2(\Psi_B(k)+\dfrac{\gamma_0}{k})$
we obtain (\ref{S1Traceformula4}).

In order to get (\ref{S1Traceformula1}), we define the
function $f_1=ik(\Psi_B(k)+\dfrac{\gamma_0}{k})$. Due to Theorem 1.1, this
function $f_1$ satisfies the condition (\ref{S2Condtiononf}) and we have
$$
\int_{\bf R}{\rm Im}\, f_1(t)dt=\int_{\bf R} t\left(\phi_B(t)+\frac{\gamma_0}{t}\right)dt=0.
$$
Substituting $\phi_B=\phi-B_I$ and $\gamma_0=\alpha_0-\beta_0$ into last integral, we obtain
$$
\int_{\bf R}t\left(\phi_B(t)+\frac{\gamma_0}{t}\right)dt=\int_{\bf R}
t\left(\phi(t)+\frac{\alpha_0}{t}\right)dt-J, \quad J=\int_{\bf R}
t\left(B_I(t)+\frac{\beta_0}{t}\right)dt.
$$
The integration by parts yields
$$
J=-\frac{1}{2}\int_{\bf R}t^{2}\left(B_I'(t)-\frac{\beta_0}{t^2}\right)dt=
-\sum_{j=0}^{N}\int_{\bf R}
t^{2}|k_j|\left(\frac{1}{t^2+|k_j|^2}-\frac{1}{t^2}\right)dt= \pi
\sum_{j=0}^{N}|k_j|^2,
$$
which implies (\ref{S1Traceformula1}).

Applying similar arguments to the function $ik^3(\Psi_B(k)+\dfrac{\gamma_0}{k})$
we obtain (\ref{S1Traceformula2}).

We will prove (\ref{S1Traceformula5}).  Below we will use the following simple
fact (see \cite{Koo}). Assume that the function $F\in \mathcal H({\bf C}_+\,;\,{\bf C})$ satisfies $ {\rm Im}\,
F(0)=0$ and $F\in L^2({\bf R})$. Then by the property of the Hilbert
transform, we have
$$
\frac{1}{\pi}\int_{{\bf R}}\frac{{\rm Im}\,F(t)}{t}dt= {\rm Re}\,F(0).
$$

We apply this result to the function $f=(\Psi_B(k)-\Psi_B(0)/(k)$. Due
to Lemma 2.2,  $f\in {\mathcal H}({\bf C}_+\,;\,{\bf C})$  and ${\rm I}\, f(0)=0$ since $\rho_B(t)$
is even on the real line. Then we have
$$
\frac{1}{\pi}\int_{{\bf R}}\frac{\rho_B(t)-\rho_B(0)}{t^2}dt=\phi_B'(0)= \phi'(+0)-B_I'(0)=\phi_{sc}'(+0)-\alpha_{-1}-2\sum _1^N\frac{1}{|k_j|},
$$
where we have used (\ref{S1detSkandDk}) and (\ref{S2logDBdiffatk=0}). This proves (\ref{S1Traceformula5}). \qed

\medskip

We prove  the uniqueness. Let $\mathcal H_+^2$
denote the Hardy class of functions $g$ which are analytic in ${\bf C}_+$ and
satisfy $\sup_{y>0}\int_{{\bf R}}|g(x+iy)|^2dx<\infty$.


\begin{lemma}
Let $f_j(z), j=1,2$, be such that $\log f_j\in \mathcal H_+^2$, and
satisfy the following conditions:

\noindent
(i) $f_j(z)=1+o(1)$ as ${\rm Im}\,z\to \infty , z\in {\bf C}_+$,

\noindent
 (ii) $\phi_j(t) := \arg f_j(t)$ is continuous in $t\in {\bf R}\setminus \{0\}$,
and $\phi_j(t)\to 0$ as $t\to \pm \infty$,

\noindent
(iii) $e^{-2i\phi_1(t)}=e^{-2i\phi_2(t)},$ for a.e. $t\in {\bf R}$.

Then $f_1(z)=f_2(z)$ on ${\bf C}_+$.
\end{lemma}

Proof. Recall that if $\phi_1(t)=\phi_2(t), $ for a.e. $t\in
{\bf R}$. Then $f_1(z)=f_2(z)$ on ${\bf C}_+$ (see [Koo]). Hence we have only to
show $\phi_1(t)=\phi_2(t), $  for a.e. $t \in {\bf R}$. We take an interval
$I=[m, m+1]$ for large $m>1$. Then $|\phi_j(t)|<\epsilon, t\in I, \ j=1,2$,
and using the identity   $e^{-2i\phi_1(t)}=e^{-2i\phi_2(t)}$ for a.e.
$t\in {\bf R}$, we have $\phi_1(t) = \phi_2(t)$ a.e. on $I$. Then by the well-known theorem on the boundary value of the harmonic functions, we have $\phi_1(z) = \phi_2(z)$ on ${\bf C}_+$, which implies $\phi_1(t) = \phi_2(t)$, a.e. on ${\bf R}$. \qed

\medskip
 Noting that  $\Psi_B(k) \in \mathcal H_+^2$ and using Lemma 2.6, we see that $D_B(k)$ is uniquely
determined by ${\arg}\,D_B(k)$ in the Hardy class.  In this class the following uniqueness
result for the required determinant $D_B(k)$ holds.


\begin{lemma}
 Let the potential $V$ satisfy
(\ref{S1ConditonC}). Let a function  $F$ be analytic in ${\bf C}_+$ and satisfy
the following conditions:

\noindent
(1) $F(k)=1+o(1)$ as ${\rm Im}\, k\to \infty , k\in {\bf C}_+$, and $\log F\in
\mathcal H_+^2$ for some branch of $\log$,

\noindent
(2) $\phi_F(t) := \arg F(t)$ is continuous in $t \in {\bf R}\setminus \{0\}$,
and $\phi_F(t)\to 0$ as $t\to \pm \infty$,

\noindent
(3) $e^{-2i\phi_F(t)}=e^{-2i\phi_B(t)},$ for a.e. $t\in {\bf R}$.

Then $F(k)=D_B(k)$ on ${\bf C}_+$.
\end{lemma}

\medskip
Note that if in this Lemma $\log F$ is not in the Hardy class, then we do not have uniqueness. For example, let $N=0$ and
$F(k)=D(k)e^{1/k},$ for  $k\in {\bf C}_+$. Then $F$ satisfies all
conditions in Lemma	 except $\log F\in \mathcal H_+^2$.

\medskip
 {\it Proof of Lemma 2.7}. We apply Lemma 2.6 to
 $D_B$ and $F$. Due  to Theorem 1.1, the
$-i\log D_B$ satisfies all conditions in Lemma 2.6. Then $D_B=F$.
\qed

\medskip
We prove our second result on the Dirichlet integrals.

\medskip
 {\it Proof of Theorem 1.2}.
 Using Green's formula we obtain the following equality, where $k = t + i\tau$.
$$
I_r(\Psi_B) := \iint_{k\in {\bf C}_+, |k|<r }|\Psi_B'(k)|^2dtd\tau =\int_{-r}^r
\phi_B'(t)\rho_B(t)dt +\frac{rJ'(r)}{2},
$$
for any $r>0$, where
$$
J(r)=\int_0^{\pi}|\Psi_B(r e^{i\varphi})|^2d\varphi.
$$
The asymptotics (\ref{S1ExpandlogDBk}) and (\ref{S2logDkprimeasymp}) yield $rJ'(r)=O(1/r), \
r\to\infty$ and $\phi_B', \ \rho_B' \in L^2({\bf R})$,
which shows that $I_r$ converges and we have
\begin{equation}
 \iint_{{\bf C}_+}|\Psi_B(k)|^2dtd\tau=\int_{{\bf R}} \phi_B'(t)\rho_B(t)dt.
 \label{S2iintPhi'BonC+}
\end{equation}
Consider the integral in the right-hand side of (\ref{S2iintPhi'BonC+}) in  more detail.
Using $\rho_B(t) = \rho(t), t \in {\bf R}$, and $\rho_B(t)\in L^1({\bf R})$  we have
the following decomposition

$$
\frac{1}{\pi}\int_{\bf R}\phi_B'(t)\rho_B(t)dt= m_B\gamma_0 +
\frac{1}{\pi}\int_{\bf R}(\phi_B'(t)-m_B)\rho_B(t)dt= m_B\gamma_0 - S_0.
$$
Here the integrals converge  absolutely, which implies (\ref{ddklogDBk}).

In order to show $\rho_B(0)=\rho(0)<0$ we use the identity for $\tau>0$
$$
\Psi(i\tau)-\Psi(0)=\int_0^{i\tau}2iz\psi (z)dz,\quad \psi (z)= {\rm Tr}\,
R_0(z)VR(z)VR_0(z)>0, \quad  z\in i{\bf R}_+.
$$
Hence $-\rho(0)=\int_0^{\infty}2y\psi (iy)dy>0,$ where the integral
converges, since $R(-\tau^2)=O(1/\tau^2)$ as $\tau \to \infty$. \qed

\section {Resonances}
\setcounter{equation}{0}

We consider first the resonances in the case when $V$ satisfies (\ref{S1ConditonC}) and $D(k)$ is entire.

\medskip
\noindent
{\it Proof of Theorem 1.3.}  Recall that $D= B D_B$. Let $D$ be entire. Then using (\ref{S1detSkandDk}), we have
the equality (\ref{S1D0-kD0k}) and then $B^2(k)\det S(k)$ has an
analytic extension from $\bf R$ into whole $\bf C_+$.

Conversely, let $B^2(k)\det S(k)$ have an analytic extension from $\bf R$ into $\bf C_+$. Then using the formula (\ref{S1detSkandDk}), we deduce that the function $D_B(-k)$ has an
analytic extension from $\bf R$ into $\bf C_+$, since $D_B(k)$ is analytic in
$\bf C_+$. Then $D$ is entire.

Now  we assume that $D$ is entire.

\medskip
\noindent
(i) The equality (\ref{S1detSkandDk}) and  $D= B D_B$ yield (\ref{S1D0-kD0k}). Then the function $\det S(k), k\in \bf R$,  has a
meromorphic extension from $\bf R$ into the whole $\bf C$.

\medskip
\noindent
(ii) Let us consider the case that the order $p=3$, the proof for the other cases is the same. The function $D$ has the form
$$
D(k)=e^{\mathcal P(k)}y(k),  \quad y(k)=\prod_1^M(k- k_n), \ \
\quad \mathcal P(k)= c_3 k^3+ c_2 k^2+ c_1 k+ C.
$$
Then the function $\Psi(k)=-i\log D(k)$ has the asymptotics
$$
\begin{aligned}
\Psi(k)=-i\mathcal P(k)-i\log y(k)=-i\mathcal P(k)-iM \log k-i\sum_{n=1}^M\log\left(1-\frac{k_n}{k}\right)\\
=-i\mathcal P(k)-iM \log k+\frac{i}{k}\sum_{n=1}^M k_n+O(1/k^2), \quad k\to \infty
\end{aligned}
$$
Then $M=0$ and $\alpha_0=0$. Hence $V=0$.

\medskip
\noindent
(iii) Using formula(\ref{S1D0-kD0k}), we obtain the statement (3).

\medskip
\noindent
 (iv) Let $Q_0(\cdot ):\bf C_+ \to \bf B_2$ have an analytic continuation
into the whole plane. Then the equality (\ref{S2I+Q0I+Q})
gives  a meromorphic continuation of $Q(k),\ k\in {\bf C}_+\setminus \sqrt{\sigma_d(H)}$, into whole $\bf C$.

If $z_0$ is a zero of $D(k) $, due to (\ref{S2I+Q0I+Q}), $z_0$ is a pole of $Q(k)$
counted with multiplicity.

 If  $z_0$ is a pole of $Q(k)$, then the equation
$(I+Q_0(z_0))f=0$ has a solution $ f\neq0$, hence $z_0$ is a zero of
$D(k)$, counted with multiplicity. \qed

\medskip
Next we consider the resonances for  compactly supported potentials $V$.

\medskip
\noindent
{\it Proof of Theorem 1.4}. The estimate in
(\ref{S1DBKinfhk}) is proved in Theorem 2.3. We prove (\ref{S1DkeC0k3}). In  \cite{Z2} there is an estimate $|B^2(k)\det
S(k)| \leq C_1e^{C_1|k|^3}$ for any $k\in \overline{\bf C_+}$ with some constant
$C_1$. Note that the proof of this last estimate is not complicated
(see also \cite{Fr}). Then using Theorem 2.3 and (\ref{S1D0-kD0k}) we
obtain (\ref{S1DkeC0k3}).

 It is well known that if $D(k)$ is entire
and has estimate (\ref{S1DkeC0k3}), then $D(k)$ has the Hadamard factorization
(\ref{S1Dkinfiniteprodut}), where $\mathcal P=Ak^3+Bk^2+Ck$. We have to determine  the constants
$A,B,C$. The function $\Psi(k)$ is odd and $\rho(k)$ is even on the real
line.
 Differentiating (\ref{S1Dkinfiniteprodut}), we get (\ref{S1D7koverD(k)}), which yields
$C=\Psi'(0)$ since $\rho'(0)=0$.  Differentiating  again (\ref{S1Dkinfiniteprodut}) we
obtain $A,B$. \qed

\medskip
Finally, we prove the result about the trace ${\rm Tr} (f(H)-f(H_0))$
in terms of resonances only.

\medskip
\noindent
{\it Proof of Theorem 1.5}.
The function $D$ is entire. Using Theorem 1.1 and (\ref{S1Kreintrace}) for each
$f \in\bf C_0^{\infty}(\bf R )$ and integrating by parts we have the following identity
$$
\begin{aligned}
{\rm Tr}\,(f(H)-f(H_0))= \sum_{j=1}^{N} f(k_j^2) -\frac{1}{\pi}\int_{\bf R_+}f(t^2)\phi_{sc}'(t)dt\\
=\sum_{j=1}^{N}f(k_j^2)-\frac{\alpha_{-1}}{\pi}
\int_{\bf R_+}f(t^2)dt-\frac{1}{\pi}\int_{\bf R_+}f(t^2)\phi'(t)dt.
\end{aligned}
$$
The expansion  (\ref{S1D7koverD(k)}) yields (\ref{S1ph7kinTheorem15}).
Using (\ref{S2logDerDk}), (\ref{S2TraceQ0'k}) and (\ref{S2I+Q0I+Q}), we have
\begin{equation}
\begin{split}
\frac{D'(k)}{D(k)}=-{\rm Tr}\, Q(k)Q_0'(k)=
-{\rm Tr} \left(Q_0'(k)-(I-Q(k))Q_0'(k)\right)\\
=-i\alpha_{-1}-2k{\rm Tr}\, (R(k^2)-R_0(k^2))
\end{split}
\nonumber
\end{equation}
and then (\ref{S1D7koverD(k)}) yields
\begin{equation}
{\rm Tr} (R(k^2)-R_0(k^2))=-\frac{D'(k)}{2kD(k)}-\frac{i\alpha_{-1}}{2k}=
-\frac{1}{2k}\left(i\alpha_{-1}+\mathcal P'(k)+\sum\frac{k^3}{k_n^3(k- k_n)}\right)
\nonumber
\end{equation}
where the series converges uniformly on any compact subset of
${\bf C}\setminus \{0, k_n, n\geq 1\}$. These two equalities imply (\ref{S1ph7kinTheorem15}).
\qed


\section {Proof of asymptotics Lemma 2.3}
\setcounter{equation}{0}

Recall that by Lemma 2.1, each function ${\rm Tr} \, Q_0^n(k), n\geq 2, $ belongs to $\mathcal H({\bf C}_+;\bf C)$
and $|{\rm Tr}\, Q_0^n(k)|\leq \|Q_0(0)\|_{2}^n$. We need the asymptotics of
${\rm Tr}\, Q_0^n(k)$ as $|k|\to \infty$. For the function $f(x),\ x\in {\bf R}^3$ and
$\omega \in S^2$, we put $\partial_{\omega} f(x)=\frac{\partial f}{\partial \omega}(x)=\omega\cdot\nabla
f(x)$.


\begin{lemma}
The following equality holds:
\begin{equation}
\frac{i}{2}{\rm Tr}\, Q_0^2(k)=-\frac{\|V\|^2}{16\pi k}-\frac{\|\nabla V\|}{3\pi (4 k)^3}
+\frac{1}{32k^4i}\int_0^{\infty} e^{2ikt}\frac{F(t)}{(1+t)^{3+\epsilon}}dt,\quad k\in {\bf C_+}.
\label{S44.1}
\end{equation}
for some $ F\in L^{\infty}(\bf R_+)$.
\end{lemma}

Proof. Let  $f(k)=\frac{i}{2}{\rm Tr}\, Q_0^2(k)$. Then we have
\begin{equation}
\begin{split}
f(k) &=\dfrac{i}{2}\iint_{{\bf R}^6}V(x)R_0(x-y,k)V(y)R_0(y-x,k)dxdy\\
&=
\frac{i}{2}\iint_{{\bf R}^6}\frac{V(x)V(y)e^{2ik|x-y|}}{(4\pi )^2 |x-y|^2}dxdy.
\end{split}
\nonumber
\end{equation}
If we set $u=x-y, v=x$ and $u=t\omega, t=|u|>0$, then we obtain
\begin{equation}
f(k)=\frac{i}{2}\int_0^{\infty} e^{2ikt}g(t)dt,\quad
g(t)=\int_{|\omega |=1}d\omega \int_{{\bf R}^3}V(x+t\omega )V(x)\frac{dx}{(4\pi )^2}.
\label{S44.2}
\end{equation}
Consider the function $g$. Using (\ref{S1ConditonC}), we deduce that
$g\in C(\bf R)$ and  $g$ is even on $\bf R$.
The derivatives of $g$ have the forms
\begin{equation}
g'(t)=\int_{S^2}d\omega\int_{{\bf R}^3}V(x)\partial_{\omega} V(x+t\omega )\frac{dx}{(4\pi )^2}
=-\int_{S^2}d\omega\int_{{\bf R}^3}V(x+t\omega )\partial_{\omega}V(x)\frac{dx}{(4\pi )^2},
\label{S44.3}
\end{equation}
\begin{equation}
g''(t)=-\int_{S^2}d\omega\int_{{\bf R}^3}\partial_{\omega}V(x+t\omega )
\partial_{\omega} V(x)\frac{dx}{(4\pi )^2}=\int_{S^2}d\omega\int_{{\bf R}^3}V(x+t\omega )
\partial_{\omega}^2 V(x)\frac{dx}{(4\pi )^2},
\label{S44.4}
\end{equation}
\begin{equation}
g'''(t)=\int_{S^2}d\omega\int_{{\bf R}^3}\partial_{\omega}V(x+t\omega )
\partial_{\omega}^2 V(x)\frac{dx}{(4\pi )^2},
\label{S44.5}
\end{equation}
where we used  integration by parts. Substituting the following estimate
$$
 \frac{1}{(1+|x|)(1+|y|)}\leq \left( \frac{1}{1+|x|}+\frac{1}{1+|y|}\right)
 \frac{1}{1+|x-y|}
$$
into (\ref{S44.2})-(\ref{S44.5}) and using (\ref{S1ConditonC}), we deduce that $g^{(n)}\in
C({\bf R}),  n=0,1,\cdots,4$, and $g^{(n)}(t)=O(t^{-3-\epsilon})$ as $t\to \pm
\infty$. We have the following equalities
\begin{equation}
\begin{split}
g(0) & =\int_{S^2}d\omega\int_{{\bf R}^3}V^2(x)\frac{dx}{(4\pi )^2} =
\frac{1}{4\pi} \int_{{\bf R}^3}V^2(x)dx,\\
g'(0)& =\int_{S^2}d\omega\int_{{\bf R}^3}V(x)\omega\cdot\nabla V(x)\frac{dx}{(4\pi )^2}
=0,\\
g''(0)&=-\int_{S^2}d\omega\int_{{\bf R}^3}(\omega\cdot\nabla V(x))^2\frac{dx}{(4\pi )^2}=
-\frac{\|\nabla V\|^2}{12\pi},
\end{split}
\nonumber
\end{equation}
since for any constant vector $\gamma \in {\bf R}^3$ we have the following
identities
$$
\int_{S^2}\gamma\cdot\omega d\omega=0,\quad
 \int_{S^2}(\gamma\cdot\omega)^2d\omega= |\gamma|^2\,2\pi \int_0^\pi \!\!\! \cos^2 \theta \sin \theta
d\theta=|\gamma|^2\frac{4\pi}{3}.
$$
Then by integration by parts, we have
$$
f(k)=-\frac{g(0)}{4k}+\frac{g''(0)}{16k^3}+\frac{1}{32k^4 i}\int_0^{\infty} e^{2ikt}g''''(t)dt,\\
$$
$$
 F(t)=(1+t)^{3+\epsilon}g''''(t)\in L^{\infty}({\bf R}_+).
$$
which implies  (\ref{S44.1}).  \qed

\medskip

We consider the function ${\rm Tr} \,Q_0^3(k)$.


\begin{lemma}
Let $\varphi(k)=
-\dfrac{i}{3}{\rm Tr}\, Q_0^3(k),\ k\in \overline{\bf C_+}$. Then we have the following asymptotic expansion:
$$
\varphi(k)=-\frac{\alpha_1^0}{k^3}+O\left(\frac{1}{k^{4}}\right),
\quad
\varphi'(k)=\frac{3\gamma_1^0}{k^4}+O\left(\frac{1}{k^{4}}\right),
$$
$$
{\rm Tr}\, Q_0^n(k)=O(k^{-4}), \quad {\rm Tr}
\left(Q_0^{n-1}(k)Q_0'(k)\right)=O(k^{-4}), \quad  n=4,5,
$$
as $|k|\to \infty,\ k\in {\bf C}_+$, where $\displaystyle{\gamma_1^0=\dfrac{2}{3(4\pi)^3}\int_{{\bf R}^3}V^3(x)dx}$.
\end{lemma}

Proof.  We have for $k \in {\bf C}_+$
\begin{equation}
\begin{split}
\varphi(k)& =-\frac{i}{3}\int_{{\bf R}^6}V(x)\frac{e^{ik|x-y|}}{4\pi |x-y|}V(y)
\frac{e^{ik|y-z|}}{\pi |y-z|}
V(z)\frac{e^{ik|z-x|}}{4\pi |z-x|}dxdydz\\
&=-\frac{i}{3}\int_{{\bf R}^9}V(z)V(z+v)V(z+u+v) \frac{e^{ik(|u|+|v|+|u+v|)}}{(4\pi
)^3|u||v||u+v|}dudvdz \\
&=-i\int_{{\bf R}^6}\frac{e^{ik(|u|+|v|+|u+v|)}}{|u||v||u+v|}f(u,v)dudv,\\
f(u,v)&=\frac{1}{3(4\pi )^{3}}\int_{{\bf R}^3}V(z)V(z+v)V(z+u+v)dz.
\end{split}
\nonumber
\end{equation}
where we used $u=x-y, v=y-z$ and $y=v+z, x=u+v+z$.
Using the new variables
$$
u=tc\omega, \quad v=ts\nu,\quad c=\cos \psi, \quad s=\sin \psi,\quad t=\sqrt {u^2+v^2}>0,\quad
\omega, \nu \in S^2,
$$
we rewrite the integral into the form
\begin{equation}
\varphi(k)=-i\int_{S^2\times S^2}d\omega d\nu
\int_0^{\pi/2}\frac{cs}{\phi} d\psi \int_0^{\infty} t^2e^{iktG} g(t,\eta )dt,\quad g(t,\eta)=f(tc\omega,ts\nu),
\label{S44.11}
\end{equation}
where
$$
\eta=(\psi, \omega, \nu)\in (0,\pi/2)\times  S^2\times S^2
,\quad \phi=|1+2cs\omega\cdot\nu|^{1\/2},\quad G=c+s+\phi.
$$
By the same arguments as in the proof of  Lemma 4.1, we deduce that each
$g^{(n)}\in C({\bf R}), \, n=0,1,\cdots,4$, and $g^{(n)}(k,\eta)=O(t^{-2(3+\epsilon)})$
as $ t\to \infty$. Moreover, we have
$$
g(0,\eta )=\int_{{\bf R}^3}\frac{V^3(x)}{3(4\pi )^3}dx,\quad
\int_{S^2\times S^2} d\omega d\nu \int_0^{\pi/2}
\frac{2cs}{\phi G^3} d\psi=2\pi ^2,
$$
$$ \int_{S^2\times S^2}
d\omega d\nu \int_0^{\pi/2} \frac{2cs}{\phi G^4}\frac{\partial g(0,\eta)}{\partial t}
d\psi=0.
$$
By integrating by parts we have
\begin{equation}
\begin{split}
& -i\int_0^{\infty} t^2e^{iktG} g(t,\eta)dt=
\int_0^{\infty} e^{iktG} \left(\frac{t^2 g(t,\eta)}{kG }\right)'dt \\
& =-\int_0^{\infty} e^{iktG} \left(\frac{t^2 g(t,\eta)}{i(kG)^2 }\right)''dt
=-\frac{2g(0,\eta)}{(kG)^3 } -\int_0^{\infty}\frac{e^{iktG}}{(kG)^3 } \left(t^2
g(t,\eta)\right)'''dt\\
& =-\frac{2g(0,\eta)}{(kG)^3 }-\frac{2}{i(kG)^4 }\frac{\partial g(0,\eta)}{\partial t}+
\int_0^{\infty} \frac{e^{iktG}}{i(kG)^4 } \left(t^2 q(t,\eta)\right)''''dt.
\end{split}
\nonumber
\end{equation}
Substituting the last integrals into (\ref{S1ConditonC}), we obtain
$$
\varphi(k)=-\frac{\alpha_1^0}{k^3} +\frac{1}{ik^4}\int_{S^2\times
S^2} d\omega d\nu \int_0^{\pi/2}\frac{cs}{\phi G^4} d\psi
\int_0^{\infty} e^{iktG} \left(t^2 g(t,\eta)\right)''''dt,
$$
which yields the first asymptotics in (\ref{S44.11}).  The proof of other
 asymptotics is similar.
 \qed

\medskip
\noindent
{\bf Acknowledgments}\ Many parts of this paper were written when E. K. was staying in
Tsukuba University, Japan. He is grateful to the Mathematical Institute for their hospitality.
This work was supported by the
Ministry of education and science of the Russian Federation, state
contract 14.740.11.0581.

\end{document}